\documentclass[12pt, reqno,fleqn]{amsart}

\usepackage{amsmath}
\usepackage{amssymb}
\usepackage{amsfonts}
\usepackage{graphicx}
\usepackage{amsthm}
\usepackage{enumerate}
\usepackage{lscape}
\usepackage{dsfont}
\usepackage{color}
\usepackage{mathtools}

\usepackage{setspace}
\newcommand{\R}{\mathds{R}}                   
\newcommand{\z}{\mathds{Z}}

\newcommand{\CP}{\mathds{C}\mathrm{P}}

\newcommand{\Ric}{\mathrm{Ric}}

\newcommand{\K}{K\"{a}hler}

\newtheorem{theor}{Theorem}

\begin{document}

\title[On the Gromov width of homogeneous  K\"{A}HLER manifold]{On the Gromov width of homogeneous K\"{A}HLER manifolds}

\author{Andrea Loi}
\address{(Andrea Loi) Dipartimento di Matematica \\
         Universit\`a di Cagliari (Italy)}
         \email{loi@unica.it}

\author{Fabio Zuddas}
\address{(Fabio Zuddas) Dipartimento di Matematica e Informatica \\
          Via delle Scienze 206 \\
         Udine (Italy)}
\email{fabio.zuddas@uniud.it}

\thanks{
The  authors were  supported by Prin 2010/11 -- Variet\`a reali e complesse: geometria, topologia e analisi armonica -- Italy. The first author was supported  by INdAM--GNSAGA - Gruppo Nazionale per le Strutture Algebriche, Geometriche e le loro Applicazioni.}

\begin{abstract}
We compute  the Gromov width of  homogeneous \K\ manifolds with second Betti number equal to  one. Our result is based on the recent preprint \cite{kk} and on the upper bound of the Gromov width for such manifolds obtained  in \cite{lmz2}.
\end{abstract}
 
\maketitle


\section{Introduction}
The Gromov width \cite{GROMOV85} of a $2n$-dimensional symplectic manifold $(M, \omega)$ is defined as
\begin{equation}\label{gromovwidth}
c_G(M, \omega)= \sup \{\pi r^2 \ |\ B^{2n}(r)\    \mbox{symplectically embeds into}  \  (M, \omega)\},
\end{equation}
where 
\begin{equation}\label{ball}
B^{2n}(r)=\{(x, y)\in\R^{2n}\  |\  \sum_{j=1}^nx_j^2+y_j^2<r^2 \}
\end{equation}
is the open ball of radius $r$ endowed with the standard symplectic form $\omega_0=\sum_{j=1}^n dx_j\wedge dy_j$ of $\R^{2n}$.
By Darboux's theorem $c_G(M, \omega)$ is a positive number.
Computations and estimates of the Gromov width for various examples have been obtained by several authors (see, e.g. \cite{lmz2} and references therein).
The main result of this paper is  the following theorem proved in the next section.

\begin{theor}\label{mainnew}
Let $(M, \omega)$ be a compact homogeneous \K\ manifold such that $b_2(M)=1$
and $\omega$ is normalized  so that $\omega(A)=\int_A\omega =\pi$ for the generator  $A\in H_2(M, \z)$. Then 
\begin{equation}\label{main}
c_G (M, \omega)= \pi.
\end{equation}
\end{theor}

The class of manifolds in Theorem \ref{mainnew} includes all  
Hermitian symmetric space of compact type whose Gromov width has been computed in 
\cite{lmz1}.  We do not know if the assumption on the second Betti number can be dropped.

\section{Proof of Theorem \ref{mainnew}}
The proof of Theorem \ref{mainnew} is mainly based  on the lower bound recently obtained by K. Kaveh
\cite{kk}:

\vskip 0.3cm
\noindent
{\bf Theorem A}
{\em Let $X$ be a smooth complex projective variety  
embedded in a complex projective space $\CP^N$. 
Then 
\begin{equation}\label{in1}
c_G (X, \omega_{FS})\geq 1,
\end{equation}}
where
$\omega_{FS}$ denotes  the restriction to $X$ of  the Fubini--Study \K\ form of $\CP^N$.

\vskip 0.3cm 

\noindent
{\bf Proof of Theorem \ref{mainnew}}
The upper bound $c_G(M, \omega)\leq\pi$ is Theorem 1 in \cite{lmz2}.
In order to obtain the lower bound $c_G(M, \omega)\geq\pi$, consider the integral \K\ form
$\hat\omega =\frac{\omega}{\pi}$ on $M$. Let $(L, h)$ be the holomorphic hermitian  line bundle on $M$ such that 
$\Ric (h)=\hat\omega$, where $\Ric(h)$ is the $2$-form on $M$ defined by
$\Ric (h)=-\frac{i}{2\pi}\partial\bar\partial\log h(\sigma, \sigma)$, for a local  trivializing holomorphic section $\sigma$ of $L$. Let $s_0, \dots ,s_N$ be an orthonormal basis for the space of global holomorphic sections  $H^0(L)$ of $L$ equipped with the  $L^2$-scalar product
 $\langle\cdot, \cdot\rangle$ given by:
$$\langle s, t\rangle =\int_Mh(s, t)\frac{\hat\omega^n}{n!},\ \  s, t\in H^0(L).$$
Then, it is not hard to see (see, e.g. \cite{arezzoloi}), due to the homogeneity and simply connectedness of $M$, that  the Kodaira map $k :M\rightarrow \CP^N, x\mapsto [s_0(x):\cdots :s_N(x)]$  is a \K\ immersion, i.e. $k^*\omega_{FS}=\hat\omega$.
Moreover, in \cite[Theorem 3]{asian} is proved that such a map is  injective, and hence 
$(M, \hat\omega)$ is symplectomorphic to $(k(M), \omega_{FS})$.
By Theorem A and by the conformality  of the Gromov width one gets
$$c_G(M, \omega)=\pi c_G(M, \hat\omega)=\pi c_G(k(M), \omega_{FS})\geq \pi$$
and the theorem is proved. 
$\hfill \Box$

\end{document}